\documentclass[12pt]{amsart}
\usepackage{amssymb}

\newcommand{\R}{\mathbb R}
\newcommand{\wh}{\widehat}
\renewcommand{\le}{\leqslant}
\renewcommand{\ge}{\geqslant}

\newtheorem*{theorem}{Theorem}
\newtheorem*{remark}{Remark}

\begin{document}

\title[]{\sc          A function with support of finite measure
                     and ``small'' spectrum}

\author{Fedor Nazarov, Alexander Olevskii}

\begin{abstract}   We construct a function on $\mathbb R$ supported 
            on a set of finite measure whose spectrum has
            density zero. 
\end{abstract}
\maketitle

\section{The result}

             Let $F$ be a function in $L^2(\mathbb R)$.  We say that it is
             supported on $S$ if 
                      $$ F= 0  \text{ almost everywhere on } \R\setminus S\,.$$
              Suppose the set $S\subset\R$ is of finite Lebesgue measure. Then
              the Fourier transform $\wh F$ of $F$ is a continuous function, so 
              the spectrum of $F$ is naturally defined as the closure 
              of the set where $\wh F$ takes non-zero values.
         
             According to the uncertianty principle, the support and 
             the spectrum of a (non-trivial) function $F$ cannot be
             both ``small sets'' .
             This principle has various versions (see e.g. \cite{HAJ94}). 

             In particular, the classic uniqueness theorem for 
             analytic functions implies that if $F$ is supported
             on an interval and it has a ``spectral gap'' 
             (that is, $\wh F = 0$ on an interval ) then $F=0$.

            Another important result says that if the support $S$ and 
            the spectrum $Q$ of $F$ are both of finite measure then $F=0$ 
            \cite{Ben74/85}, \cite{AB77}.

            On the other hand, $F$ may have a support of finite measure 
            and a spectral gap; see \cite{Kr82}, where such an example was
            constructed with $F= 1_S$.

             Answering a question posed by Benedicks, Kargaev 
            and Volberg \cite{KV92} constructed an example of a function
            $F$ such that
           $$ |S| < +\infty,   |\R \setminus Q| = +\infty $$
            (here and below by $|A|$ we denote the Lebesgue measure of 
             the set $A$).

            The goal of this note is to prove the  following

\begin{theorem}    
There is a function $F\in L^2(\R)$ supported by 
            a set $S$ of finite measure, such that 
$$
                       |Q\cap (-R,R)| = o(R) \text{ as }R\to \infty\,.
$$
            In addition, $F$ can be chosen as the indicator function of $S$.
\end{theorem}

            The proof below is based on a simple construction, 
            completely different from the ones in the cited papers.

\section{Proof}

\subsection{}    Take a Schwartz function $F_0$ such that
$$
                         0 \le F_0 (t) \le 1\qquad  (t\in \R)
$$
             and its Fourier transform $\wh F_0$ is positive on 
             $(-1,1)$  and vanishes outside that interval.              
             Define a sequence of functions $F_n$  recursively by
$$
                         F_n := F_{n-1} + G_n  \qquad (n=1,2,...)\,,                                   
$$
             where
\begin{equation}
\label{eq1}
             G_n(t) :=  F_{n-1}(t) [1-F_{n-1}(t)] \cos k_n t
\end{equation}
                     
            We are going to prove that if the numbers $k_n$ grow 
            sufficiently fast, then the sequence $F_n$ converges to
            a function $F$ satisfying the requirements of 
            the theorem.
          
\subsection{}   Clearly, $F_n$ and $G_n$ are Schwartz functions.

             A simple induction shows that for every $t\in \R$, we have
$$
              |G_n (t)| \le \max \{ F_{n-1} (t), 1- F_{n-1} (t)\}                  
$$
              and
$$
                    0 \le F_n (t) \le 1 \, .                                        
$$                         

              The Fourier transforms of $F_{n-1} [1-F_{n-1}]$, $F_{n-1}^2 [1-F_{n-1}]$, and $F_{n-1}^2 [1-F_{n-1}]^2$ vanish
             outside a compact interval, so for each $n\ge 1$, we have:
$$  
                 \int_{\R} G_n=\int_{\R} F_{n-1}G_n =0 
$$
and
$$
\int_{\R} G_n^2=\frac 12\int_{\R} F_{n-1}^2 [1-F_{n-1}]^2\,,
$$
              provided that $k_n$ is chosen sufficiently large.                                     
              It follows that
\begin{equation}
\label{eq2}
             \int_{\R} F_n = \int_{\R} F_0 =: C                                        
\end{equation}
             and, thereby,
$$
                    I_n:= \int_{\R} F_n(1-F_n)\le C                                 
$$
               (here, as usual, by $C$ we denote a positive constant that may 
              vary from line to line).

               Observe also that
$$
               I_n  =\int_{\R} [F_{n-1} + G_n] [1-F_{n-1} -G_n] = I_{n-1} - \int_{\R} G_n^2,
$$
               which implies that      
$$
                 \sum_{n\in[1,N]} \int_{\R} G_n^2 \leq  I_0 - I_N \le C,
$$
               and so                          
\begin{equation}
\label{eq3}
                         \sum_n  \int_{\R} G_n^2 \le C     
\end{equation}
           
\subsection{}   Define the sequence  $Q_n$ of intervals on (another copy of)
                $\R$ recursively as follows: 
\begin{align*}
                 Q_0&:=[-1,1],    
\\
                 Q_n&:= \operatorname{conv}(Q_{n-1} \cup [k_n + 2 Q_{n-1} ] \cup 
                         [-k_n + 2 Q_{n-1}])      
\end{align*}
(here $\operatorname{conv}E$ denotes the convex hull of a set $E\subset\R$).
                Clearly,  for every n,
\begin{align*}
                        \operatorname{spec} F_{n-1}&\subset Q_{n-1}\, ;
\\
                        \operatorname{spec} G_n &\subset   [k_n +2 Q_{n-1}] \cup  [-k_n+ 2 Q_{n-1}].
\end{align*}
                  Set  $Q:= Q_0\cup \bigcup_n ([k_n +2 Q_{n-1}] \cup  [-k_n+ 2 Q_{n-1}])$.
 
                 Choosing $k_n$ growing sufficiently fast we can ensure that the
                 spectra of $G_n$ are parewise disjoint and
\begin{equation}
\label{eq4}
                            |Q \cap (-R,R) | = o(R)\text{ as } R\to\infty\,.
\end{equation}                   
        
\subsection{}  Consider the series $F_0 + G_1 + G_2 +\dots$
                 Since the spectra of the terms are pairwise disjoint,
                 this series is orthogonal in $L^2(\R)$. 
                 Then (\ref{eq3}) implies that it converges in $L^2(\R)$
                 to some non-trivial function $F$.
                 The partial sums of this series are $F_n$.  
                 Take a subsequence $F_{n_\ell}$ such that  
$$                
                      F_{n_\ell} \to F  \text{ almost everywhere on } \R \text{ as } \ell\to\infty\,. 
$$
                 Recall that all $F_n$ are non-negative functions, so (\ref{eq2})
                 implies that
\begin{equation}
\label{eq5}
                        F\ge 0\text{ almost everywhere and }   
                        \int_{\R} F < \infty.                
\end{equation}        
                It follows from (\ref{eq1}) and (\ref{eq3}) that  
$$
                     \sum_n \int_{\R} [F_n(1-F_n)]^2  = 2 \sum_n \int_{\R} G_n ^2
                     < +\infty\,, 
$$
                so we must have
$$                
                F(1-F)=\lim_{\ell\to\infty} F_{n_\ell}(1-F_{n_\ell})  = 0\text{ almost everywhere,}
$$
                which implies that $F$ is the indicator-function
                of a set $S$.
According to (\ref{eq5}), this set has finite measure.
                Clearly the spectrum of $F$ is a subset of $Q$.            
                Due to (\ref{eq4}) it has density zero. This finishes 
                the proof.    
\begin{remark}
              Consider the function
$$
                  h(R):= |Q \cap (-R,R)|.
$$
              In the conditions of the Theorem, it can not be bounded.
              However the proof above shows that it may increase
              arbitrarily slowly.   
It remains an open question, however, if $Q$ can have {\em uniform} density
$0$, i.e., if it is possible that 
$$
\lim_{R\to\infty}\sup_{x\in\R}\frac 1{2R}|Q \cap (x-R,x+R)|=0\,.
$$
\end{remark}

\end{document}